\documentclass[12pt]{article}
\usepackage[english]{babel}
\usepackage{amsfonts}
\usepackage{amsmath}
\usepackage{amssymb}
\usepackage{mathrsfs}
\usepackage{latexsym}
\usepackage{multicol}
\usepackage{fancybox} 
\usepackage{color}
\usepackage{hyperref}
\usepackage{graphicx} 
\usepackage{caption}

\usepackage{upgreek}
\def\be{\begin{equation}}
\def\ee{\end{equation}}
\def\tl{\tilde} 
 
\def\gm{\gamma} 
\def\lm{\lambda}

\def\d'{``}

\newtheorem{thm}{Theorem}[section]

\newtheorem{propn}[thm]{Proposition}

\newtheorem{rem}[thm]{Remark}

\newtheorem{conje}[thm]{Conjecture}

\newtheorem{lem}[thm]{Lemma}

\newenvironment{prf}{\trivlist \item [\hskip
\labelsep {\bf Proof:}]\ignorespaces}{$\Box$ \endtrivlist}

\def\be{\begin{equation}}
\def\ee{\end{equation}}
\def\bea{\begin{eqnarray}}
\def\eea{\end{eqnarray}}

\def\i'{\textrm{i}}

\def\d'{``}
\begin{document}

%\begin{center}
%\Large{\bf{Properties of the series solution for Painlev\'e I}}
%\end{center}
\title{Properties of the series solution for Painlev\'e I} 
\author{A.N.W. Hone\thanks{SMSAS, University of Kent, Canterbury, U.K.}, 
O. Ragnisco\thanks{Dipartimento di Fisica, Universit\`a Roma Tre, Roma, Italy.} 
and F. Zullo$^{*\dagger}$  }

%\begin{center}
%{ {A.N.W. Hone\thanks{SMSAS, University of Kent, Canterbury, U.K.}, 
%O. Ragnisco\thanks{Dipartimento di Fisica, Universit\`a Roma Tre, Roma, Italy.
%and F. Zullo  }}
%{ }
%\end{center}

%\medskip
%\medskip

%\bigskip\bigskip

%\noindent

%\noindent
%KEYWORDS: Laurent expansion, Painlev\'e equation, tau-function, sigma function.  
\maketitle 

\begin{abstract}
\noindent
We present some observations on the asymptotic behaviour of the coefficients in the Laurent 
series expansion of %certain 
solutions of the first Painlev\'e equation.  For the general solution, explicit recursive formulae for the 
Taylor expansion of the tau-function around a zero are given, which are natural extensions of analogous formulae 
for the elliptic sigma function, as given by Weierstrass. 
Numerical and exact results on the symmetric solution which is singular at the origin %pentagonal symmetry 
are also presented. 
\end{abstract} 

\section{Introduction} \label{intro}
\vspace*{4mm} 

The first Painlev\'e equation ($P_I$) is usually written in the canonical form 
\be \label{p1orig} 
\frac{d^2 u}{dz^2} = 6u^2 + z. 
\ee 
It is known that all  solutions of (\ref{p1orig}) are non-classical transcendental functions, in the sense that they cannot 
be written in terms of algebraic or elliptic functions of $z$, nor in terms of classical special functions defined by 
linear differential equations. The $P_I$ equation was discovered in Painlev\'e's classification of second 
order differential equations  %(of a specific form) 
whose solutions have the property that %of being single-valued in the neighbourhood of 
all movable singularities are poles. 
For the equation (\ref{p1orig}), 
it is further known that all solutions of this equation are meromorphic, with an infinite number of poles 
in the complex plane. 
Although Painlev\'e's original approach to proving meromorphicity, as outlined in chapter XIV of Ince's book \cite{ince}, 
appears to have some gaps from the modern point of view, these issues have subsequently been resolved 
(see  \cite{hl}), and other rigorous proofs are also available, either indirectly by reformulating an 
associated isomonodromy problem as a regular Riemann-Hilbert problem \cite{fokas},  
%with the most concise treatment being that of Steinmetz 
or more directly using differential inequalities \cite{steinmetz}.  

The papers in this special issue are written in honour of Okamoto's insights 
concerning the space of initial values for the Painlev\'e equations. The important 
point made by Okamoto is that, in order to consider the solutions of a system with poles, such as (\ref{p1orig}), 
one should enlarge the phase space: as well as regular initial conditions $(u(z_0),u'(z_0))\in\mathbb{C}^2$, 
it is necessary to consider the case where $u$ (hence also $u'$) has a pole at $z=z_0$; this requires adding points 
at infinity to the phase space, and leads to a sequence of blowups around singular points \cite{joshi_hans}. 

Boutroux showed that, for %suitable $z$ with 
large $|z|$, 
with the scaled variables $U=z^{-1/2}u$ and $Z=\frac{4}{5}z^{5/4}$, 
the solution of (\ref{p1orig}) behaves asymptotically 
like an elliptic function, $U\sim \wp$, where $\wp$ is the Weierstrass  function\footnote
{Precisely, the asymptotics on fixed rays corresponds to parameters 
$g_2=12$ and $g_3= g_3(\arg z)$.}, %where $g_3$ is different on different rays.}, 
which satisfies the second order differential equation 
\be\label{wp} 
\wp '' = 6\wp^2 -\frac{1}{2}g_2. 
\ee 
%(To be precise, the asymptotics on fixed rays corresponds to parameters 
%$g_2=12$ and $g_3\equiv g_3(\arg z)$, where $g_3$ is different on different rays.)  
The purpose of this article is to present some remarks about the direct comparison between the solutions of 
the $P_I$ equation and the Weierstrass elliptic functions, at the level of exact series expansions rather 
than asymptotics. 

After introducing our conventions of scaling and notation in the next section, in section 3 we  proceed to compare 
the coefficients in the Laurent series for the $\wp$ function, which are written in terms of modular forms, with 
the corresponding expansion of solutions of (\ref{p1orig}) around a pole. This produces a natural analogy between 
the lemniscatic and equianharmonic elliptic functions, which have the symmetries of the square and the hexagon, respectively,
and a special solution of $P_I$ which is singular at the origin and has pentagonal symmetry. In section 4 we apply the same idea to the tau-function 
of $P_I$, and show how both the Hirota bilinear form, and an associated equation of degree four, lead 
to recursive formulae for the Taylor coefficients in the expansion around a zero. This both generalizes the results of 
Eilbeck and Enolskii  on the Weierstrass sigma function \cite{ee}, and at the same time extends some classical 
expansion formulae due to Weierstrass \cite{weierstrass} to the case of $P_I$. As well as being of theoretical  
interest, we show that these formulae are also useful for doing numerical computations.

\section{The equation with parameters}\label{sec1}

The programme of classification %that was  
initiated by Painlev\'e began by considering the necessary  
conditions for an equation such as (\ref{p1orig}) to have only poles as movable singularities. This requires 
some local analysis to identify the form that algebraic singularities can take at leading order, followed by 
a power series expansion around any such singularity to check whether any secular terms can appear; 
nowadays this procedure is often referred to as Painlev\'e analysis (see \cite{hone} and references). 
For the equation (\ref{p1orig}) the only algebraic singularities are double poles, 
$y\sim (z-z_0)^{-2}$, where $z_0$ is the (movable) position of the pole. This leading order behaviour 
then extends to a power series, whose resonances (the places where arbitrary coefficients can appear 
in the expansion) are found by substituting 
$
u\sim (z-z_0)^{-2} \left( 1 + \epsilon (z-z_0)^r \right) 
$ 
into the equation and comparing the leading order terms that are linear in $\epsilon$. In this case one finds 
a quadratic equation in $r$, with roots $-1,6$; the root $r=-1$ corresponds to the fact that the position of the pole at 
$z_0$ is movable, 
while  the value $r=6$ means that the coefficient of $(z-z_0)^4$ in the expansion should be arbitrary. 
Finally, one must %then 
calculate the terms of a full Laurent expansion around the pole up to 
order $(z-z_0)^4$, and  check that no additional counterterms involving $\log (z-z_0)$ are 
needed for consistency. Hence one obtains a local series representation of the general solution, 
with two free constants corresponding to the resonances.   
Assuming that such a Laurent expansion is convergent in some neighbourhood of $z_0$, 
this verifies that (\ref{p1orig}) 
has solutions that are locally meromorphic around a pole. 
Then for $P_I$, any local Laurent expansion around a pole completely fixes a global solution of the equation, 
by analytic continuation.  
However, purely local analysis is not sufficient to establish  
that all solutions are meromorphic in $\mathbb{C}$, since global estimates are necessary to ensure that 
poles cannot coalesce to form an essential singularity (see e.g. \cite{steinmetz}).

For what follows, it is worth mentioning at this stage that the details of the Painlev\'e analysis 
for the differential equation (\ref{wp}) are almost identical to those for $P_I$. However, if one 
replaces   the $z$ in (\ref{p1orig}), or the constant $-\frac{1}{2}g_2$ in (\ref{wp}), by 
an arbitrary holomorphic function $f(z)$, then %for consistency 
at the resonance $r=6$ 
one requires $f''(z_0)=0$ for all $z_0$, and hence $f$ must be a linear function,  
in order to have a consistent Laurent series. 

In the rest of the paper we consider solutions of (\ref{p1orig}) with a pole at $z_0$, so it 
is convenient to replace $z\to z+z_0$, which adds a constant to the right hand side of (\ref{p1orig}). 
Furthermore, for comparison with the equation (\ref{wp}) it is useful to rescale 
$u\to (-6\lm )^{-2/5} u$
 and 
$z\to (-6\lm )^{1/5}z$, 
so that we 
arrive at the following equation: 
%Let us consider the Painlev\'e I equation:
\be\label{PI}
\frac{d^2 u}{d z^2} = 6u^2-6\lm z-\frac{g_2}{2}. 
\ee
The above equation has two parameters, namely $\lm$ and $g_2$, 
and any solution of the  $P_I$ equation in the canonical form  (\ref{p1orig}) with a double pole at $z=z_0$ 
corresponds to a solution of (\ref{PI}) with a pole at $z=0$, for  a suitable choice of 
the constant $g_2$,  with any $\lm\neq 0$. Moreover, when $\lm =0$ the equation (\ref{PI}) reduces to (\ref{wp}). 

\subsection{Recursion for Laurent coefficients} 

To determine the local Laurent expansion of the solution of (\ref{PI}) 
with a double pole at the origin, we insert the expression 
\begin{equation}\label{series} 
u(z)=\sum_{n=0}c_n z^{n-2},  \qquad c_0=1, 
\end{equation}
into the equation, from  which
% it follows that the leading coefficient is  $$ $$ while the 
subsequent coefficients 
%(\ref{PI}), 
are determined by the recursion 
\be \label{orec}
(n+1)(n-6) c_n = 6 \sum_{j=1}^{n-1} c_j c_{n-j} -\frac{1}{2}g_2\delta_{n,4}-6\lm \delta_{n,5}, \qquad n\geq 1. 
\ee 
The above relation determines all the coefficients $c_1,\ldots ,c_5$ uniquely, but $c_6$ is arbitrary 
(corresponding to the aforementioned resonance). Once $c_6$ has been fixed, then all $c_n$ for $n\geq 7$ 
are all uniquely determined by (\ref{orec}). We find 
\be \label{inits}
c_1=c_2=c_3=0, \quad c_4=\frac{g_2}{20}, \quad c_5=\lm, \quad  c_6=\frac{g_3}{28}, \quad c_7 = 0, \quad \ldots , %\qquad 
\ee
where the value of $c_6$ is given in terms of an additional parameter $g_3$, which is introduced 
for ease of comparison with the elliptic case.   
The following result is immediate. 

\begin{lem} For $n\geq 1$ the coefficients of the expansion (\ref{series}), with $c_6=g_3/28$, are polynomials 
in $g_2$, $\lm$ and $g_3$,  
i.e. $
c_n = \mathrm{P}_n (g_2,\lm ,g_3 )$,  
where $\mathrm{P}_n$ has rational coefficients and is weighted homogeneous of total degree $n$ in these 
arguments with weights $4,5$ and $6$ respectively:  
$$ 
\mathrm{P}_n (\zeta^4 g_2,\zeta^5 \lm ,\zeta^6g_3 ) = \zeta^n \, \mathrm{P}_n (g_2,\lm ,g_3 ) 
\qquad \forall \zeta \in \mathbb{C}^* .
$$ 
\end{lem} 

The polynomials $\mathrm{P}_n$ for $n=1,\ldots,7$ can be read off from (\ref{inits}). Note that, a priori, they 
must be identically zero for $n=1,2,3,7$, since there are no non-zero homogeneous polynomials with these 
weights; but they are non-zero for all other values of $n$.  Here we list the next few non-trivial ones: 
$$ 
\mathrm{P}_8 = \frac{g_2^2}{1200}, \quad \mathrm{P}_9=\frac{g_2\lambda}{50}, \quad 
 \mathrm{P}_{10}=\frac{3}{22}\left( \lambda^2+\frac{g_2g_3}{280}\right),  
\quad \mathrm{P}_{11}=\frac{g_3\lambda}{140},
$$ % 
$$  
\mathrm{P}_{12}=\frac{1}{208}\left( \frac{g_2^3}{750}+\frac{g_3^2}{49}\right), 
\quad 
\mathrm{P}_{13}= \frac{11g_2^2\lambda}{49000}, 
\quad 
\mathrm{P}_{14}=\frac{1}{44}\left( \frac{59g_2\lambda^2}{500}+\frac{g_3g_2^2}{4200}\right).
$$
These polynomials  have some interesting properties. As will be explained in section 3, for 
$\lm =0$ the polynomials $  \mathrm{P}_{2n} (g_2,0 ,g_3 )$ are polynomials in $g_2$ and $g_3$, which 
define modular forms: up to overall scaling they give the Eisenstein series of %corresponding to 
the elliptic curve 
\be \label{curve} 
 y^2=4x^3-g_2 x-g_3 . 
\ee
In fact, viewing  each $\mathrm{P}_n$ as a polynomial in $\lm$ with coefficients in $\mathbb{Q}[g_2,g_3]$, 
each coefficient is weighted homogeneous in $g_2,g_3$ with some weight, so also defines a modular form. 
Thus we refer to the $\mathrm{P}_n$ as {\it modular polynomials}. 

Henceforth we let $u(z;g_2,\lm ,g_3 )$ denote the meromorphic function defined by the series (\ref{series}). 
Our main task in  section 3 will be to consider the asymptotic behaviour of the coefficients $c_n$ in this series 
as $n\to \infty$, and the way that this depends on the parameters $g_2,\lm , g_3$. Before doing this, 
we introduce some objects needed in section 4. 

\subsection{Hamiltonian and tau-function}

Hamiltonian forms of the Painlev\'e equations were originally found by 
Malmquist \cite{malmquist}. Some years later, Okamoto introduced the notion of 
a tau-function associated with the Hamiltonian \cite{okamoto}. 

The equation (\ref{PI}) can be written as a Hamiltonian system like so: 
\be \label{ham}  
u' = \frac{\partial h}{\partial v}, \qquad v'=-\frac{\partial h}{\partial u}, 
\qquad \mathrm{with} \qquad 
h=\frac{1}{2}v^2-2u^3+\frac{1}{2}g_2 u +6\lm z u + \frac{1}{2}g_3. 
\ee 
(The constant $g_3$ is included in the Hamiltonian $h$ for later convenience.) 
%(The inclusion of the constant $g_3$ does not affect Hamilton's equations, but it will be convenient 
%for what follows.) %$$ 
%$$ where (with t )the Hamiltonian is %\footnote{The constant $g_3$ is included for later convenience, but does not affect Hamilton's equations.} 
Taking the total derivative of $h$ with respect to $z$ we have 
\be \label{hp} 
h' = \frac{\partial h}{\partial z} = 6\lm u. 
\ee 
The tau-function associated with a solution of (\ref{PI}) 
is a function $\uptau =\uptau (z)$  such that %defined by the equation 
\be\label{tau} 
 u = -\frac{d^2}{dz^2} \log \uptau . 
\ee 
Observe that the preceding relation only defines $\uptau$ up to %the action of 
gauge transformations 
\be \label{gauge} 
\uptau \rightarrow \exp (a z +b) \, \uptau, \qquad a,b \quad \mathrm{arbitrary}, 
\ee 
and the fact that $u$ is meromorphic means that $\uptau$ is holomorphic,  with $\uptau$ having  
simple zeros wherever $u$ has double poles.   

We also note from (\ref{hp}) that, up to the addition of %fixing an integration 
a constant, we have 
\be \label{htau} 
h=-6\lm \frac{d}{dz}\log \uptau .
\ee 
The constant of integration  has been fixed so that when $\lm =0$ we have $h=0$,  
and the right hand side of (\ref{ham}) corresponds to the first order differential equation 
for the Weierstrass $\wp$ function, i.e. 
\be \label{wpeq} 
(\wp ')^2 =4\wp^3 -g_2\wp -g_3, %constant 
\ee 
for $u=\wp$ and $v = \wp '$. This choice of integration constant will also be convenient for the 
analysis of the tau-function when $\lm \neq 0$. 

\section{Asymptotics of coefficients in the Laurent series} 

We now describe the behaviour of the coefficients $c_n$, which satisfy the recursion 
\begin{equation}\label{rec1} % \left\{\begin{aligned}
c_n=\frac{6}{(n+1)(n-6)}\sum_{k=1}^{n-1}c_kc_{n-k},  \qquad  n\geq 7, 
\end{equation} 
for a suitable choice of initial values. 
To begin with we consider the classical case 
$\lm=0$, with $u(z;g_2,0,g_3 )=\wp (z; g_2,g_3)$, 
for which $c_1=c_2=c_3=c_5=0$ in (\ref{rec1}).

%Furthermore, all 
%important insight 

\subsection{The elliptic case}  %$\lm=0$.} 
\label{subsec1}  

The Weierstrass  $\wp$ function can be defined 
as the unique solution of the first order differential equation (\ref{wpeq}) having  
a double pole at $z=0$. However, in most standard treatments of 
elliptic function theory (see e.g. chapter XX in \cite{WW}),  
it is usually defined by %starting from 
its Mittag-Leffler expansion, that is 
% for $\wp(x,g_2,g_3)$ is given by:
\be\label{ML}
\wp(z; g_2,g_3) = \frac{1}{z^2}+\sum_{\Omega\in \Lambda\setminus \{ 0\} }\left(\frac{1}{(z-\Omega)^2}-\frac{1}{\Omega^2}\right) , 
\ee
where $\Lambda$ denotes the period lattice, consisting of all periods of the form $\Omega=2m\omega_1+2n\omega_2$ for 
$(m,n)\in\mathbb{Z}^2$, generated by the two independent periods 
$2\omega_1$, $2\omega_2$. The uniform convergence of the expansion (\ref{ML}) 
is guaranteed by the fact that the power 
sums $\sum_{\Omega \neq 0} \Omega^{-\mu}$ 
are absolutely convergent for all $\mu>2$, and expanding 
around $z=0$ gives the series (\ref{series}) with $c_0=1$, $c_1=c_2=0$ and  
\be \label{cexp} 
c_n=(n-1) \sum_{\Omega \neq 0} \Omega^{-n}, \qquad  n\geq 3. 
\ee
Thus it follows from 
(\ref{ML}) that $\wp$ satisfies the differential equation (\ref{wpeq}) with 
$$ 
g_2 = 60 \sum_{\Omega \neq 0} \Omega^{-4}, \qquad 
g_3 = 140 \sum_{ \Omega \neq 0} \Omega^{-6}. 
$$ 
Conversely, given a solution of (\ref{wpeq}) with parameters 
$g_2,g_3$, two independent periods are obtained from the integrals 
$ 2\omega_j = \oint_{\mathfrak{c}_j} \frac{dx}{y}$, $j=1,2$ 
around %two independent 
cycles $\mathfrak{c}_1,\mathfrak{c}_2$ that 
generate the homology of the elliptic curve (\ref{curve}). 
(This is in the generic case $g_2^3-27g_3^2 \neq 0$; otherwise 
(\ref{wpeq}) is solved in elementary functions.)  
  
%To understand the asymptotic behaviour of the coefficients $c_n$, first note 
Observe that, from  
the symmetry of the lattice $\Lambda = -\Lambda$, where the 
poles of $\wp (z)$ are situated, the odd index coefficients 
$c_{2n+1}=0$ for all $n\geq 1$ (hence $\wp$ is an even function).  
%If $2\omega_1$ is a period of minimum modulus, then we  
The even index sums (\ref{cexp}) can be rewritten as  
\be\label{cform} 
c_{2n}=\frac{2n-1}{(2\omega_1)^{2n}}\, G_{2n}(\tau ) ,
\ee 
where the function $G_{2n}(\tau)$ is the Eisenstein series, 
a modular form of weight $2n$ 
corresponding to a normalized lattice with periods 
$1,\tau$, with $\tau = \frac{\omega_2}{\omega_1}$:  
%such that 
\begin{equation}\label{eis} 
G_{2n}(\tau)=\sum_{(p,q)\in\mathbb{Z}^2\setminus \{(0,0)\}} 
\frac{1}{(p+\tau q)^{2n}}. 
\end{equation}

The periods $2\omega_{1}$, $2\omega_{2}$  %$\tau$ 
can always be ordered so that Im$\,\tau >0$, 
and by choosing  $|2\omega_1|$ to be minimal we also 
have $|\tau |\geq 1$. Then the series (\ref{series}) 
has radius of convergence  $|2\omega_1|$, and the asymptotic behaviour 
of $c_{2n}$ depends on $\tau$. 
There are two possibilities: 
\begin{itemize} 
\item ${\mathbf |\tau |>1}${\bf :} $\lim_{n\to\infty} G_{2n}(\tau )=2$; 
\item ${\mathbf |\tau |=1}${\bf :} $G_{2n}(\tau )$ has no limit as $n\to\infty$.
\end{itemize}    
Note that we can always choose $\tau$ so that it lies in the fundamental 
domain 
$\mathcal{F}=\{ |\tau|\geq 1\}\cap\{ |\mathrm{Re}\,\tau|\leq\frac{1}{2}\}$ 
(\cite{silverman}, Proposition 1.5). In the first situation, $\pm 1$ 
are the only 
normalized periods that lie on the unit circle, so all of the terms 
in the sum (\ref{eis}) vanish as $n\to \infty$ apart from the 
contribution from $(p,q)=(\pm 1 ,0)$. For the second possibility, 
$\tau =e^{\i'\theta}$ lies on the lower boundary of $\cal{F}$, 
and (by using the modular transformation $\tau\to -1/\tau$) it can 
always be assumed that $\pi /3\leq \theta \leq \pi /2$; then   
$\theta = \pi /3$, with six normalized periods on the unit circle 
(the sixth roots of unity), or otherwise there are four periods 
$\pm 1,\pm \tau$ with modulus one. The extreme cases are 
the so called \emph{equianharmonic} case ($\theta = \pi /3$), 
and the \emph{lemniscatic case} ($\theta = \pi /2$).  

\noindent 
{\bf Equianharmonic case:} 
This corresponds to $g_2=0$ in (\ref{wpeq}). 
Without loss of generality (by rescaling) 
one may set $g_3=1$ and 
find half-periods 
\be \label{omeqan}
\omega_1=\int_{\sqrt[3]{4}}^\infty\frac{dx}{\sqrt{4x^3-1}} 
=\frac{\Gamma^3\left(\frac{1}{3}\right)}{4\pi} 
, % g_3^{\frac{1}{6}}}, 
\qquad \omega_2=e^{\frac{\i'\pi}{3}}\omega_1, 
\ee
which gives $\tau=e^{\frac{\i'\pi}{3}}$. 
Then from the recurrence (\ref{rec1}) with $c_4=0$ 
it follows that $G_{2n}(e^{\frac{\i'\pi}{3}})=0$ unless 
$n =0\, \bmod \, 3$.   %either directly 
The Laurent series for the $\wp$ function %in this case 
becomes 
%(x,0,g_3)$ can be written as:
\begin{equation}\label{ex1}
\wp(z;0,1) = \frac{1}{z^2}+\sum_{n=1}^{\infty}
\frac{(6n-1)}{(2\omega_1)^{6n}}\, G_{6n}(e^{\frac{\i'\pi}{3}})\, 
z^{6n-2} 
\end{equation}
with $\omega_1$ as in (\ref{omeqan}).  
There are six normalized periods 
$e^{\frac{\i' j\pi}{3}}$, $j=0,\ldots ,5$ on the unit circle, 
%from which it follows 
so in the limit 
%Note also that the corresponding values 
of the Eisenstein series (\ref{eis}),  only six terms survive to yield 
$\lim_{n\to \infty}G_{6n}(e^{\frac{\i'\pi}{3}}) =6$.    
Some %of the 
values of these modular functions 
are given  %we report some values of the series up 
to twenty decimal places in Table \ref{tab1}. %digits of precision.  
\begin{table}
\begin{center}
\begin{tabular}{|c|c|}
\hline
$n$ & $G_{6n}(e^{\frac{\i'\pi}{3}})$\\
\hline \hline
1 & 5.86303169342540159797 \\
2 & 6.00963997169768048102\\
3 & 5.99971835637052593409\\
4 & 6.00001164757977973485\\
5 & 5.99999958743553301523\\
6 & 6.00000001557436652006\\
... & ...\\
11 & 5.99999999999999892076\\
12 & 6.00000000000000003997\\
13 & 5.99999999999999999851\\
14 & 6.00000000000000000005\\
\hline
\end{tabular}
\caption{The oscillating values of $G_{6n}(e^{\frac{\i'\pi}{3}})$.}
\label{tab1}
\end{center}
\end{table}

\noindent
{\bf Lemniscatic case: } 
This case corresponds to $g_3=0$, and (upon scaling so that 
$g_2=4$)  the half-periods are found as 
$$ 
\omega_1 = \int_{0}^{1}\frac{dt}{\sqrt{1-t^4}} 
= \frac{1}{4}B(1/4,1/2), \qquad \omega_2 = \i' \omega_1 ,  
$$  
whence $\tau =\i'$. 
The fact that $c_6=0$ in (\ref{rec1}) 
now implies that 
$G_{2n}(\i' )=0$ unless $n$ is even, so that the series (\ref{series}) 
in this case takes the form 
\be\label{ex2}
\wp(z;4,0) = \frac{1}{z^2}+\sum_{n=1}^{\infty}
\frac{4n-1}{(2\omega_1)^{4n}}\, G_{4n}(\i')\, z^{4n-2}.  
\ee
The coefficients can also be expressed in terms of the %so-called 
Hurwitz numbers $H_n$ 
\cite{Carlitz, Hurwitz, LRL} (elliptic analogues of the Bernoulli numbers) 
which are %explicitly 
given by %\cite{Hurwitz}:
%\begin{equation}
$H_n=\frac{(4n)!}{(4w)^{4n}}G_{4n}(\i')$. 
%\end{equation}
From (\ref{rec1}) 
it follows that these numbers satisfy the %showed that the numbers $H_n$ satisfy the 
recurrence found by Hurwitz:
\be\label{Hrec}
H_n=\frac{3}{(2n-3)(16n^2-1)}\sum_{k=1}^{n-1}(4k-1)(4n-4k-1){4n \choose 4k}H_kH_{n-k} . 
\ee
As is suggested by Table \ref{tab2}, the limit of the 
%Now the corresponding values of the 
Eisenstein series (\ref{eis}) in this case is 
$
\lim_{n\to \infty}G_{4n}(\i' ) =4$.  
%rapidly converge to the numerical value $4$: in  
%we report some values of the series up to twenty decimal digits of precision.  
\begin{table}
\begin{center}
\begin{tabular}{|c|c|}
\hline
$n$ & $G_{4n}(\i')$\\
\hline \hline
1 & 3.15121200215389753821 \\
2 & 4.25577303536518951844 \\
3 & 3.93884901282797037475 \\
4 & 4.01569503302502485587 \\
5 & 3.99609675317628955957 \\
6 & 4.00097680530383862810 \\
... & ...\\
11 & 3.99999904632591103400 \\
12 & 4.00000023841859318284 \\
13 & 3.99999994039535611558 \\
14 & 4.00000001490116124950 \\
%...& ...\\
\hline
\end{tabular}
\caption{The oscillating values of $G_{4n}(\i')$.}
\label{tab2}
\end{center}
\end{table}

%\end{description}

\subsection{The general case $\lm \neq 0$}\label{sec2} 

All solutions of $P_I$ are known to have order of growth 
$5/2$ \cite{steinmetz2}, and this is the same as the infimum of the values 
of $\mu$ such that the power sums over non-zero poles, 
$\sum_{\Omega\neq 0} \Omega^{-\mu}$, are absolutely convergent 
(\cite{titchmarsh}, chapter VIII).
As a consequence, the solutions of (\ref{PI}) with $\lm\neq 0$ 
admit a Mittag-Leffler expansion of the same form (\ref{ML}) 
as for $\wp(z)$: %in the elliptic case: % that is 
\be\label{MLpi}
u(z; g_2,\lm ,g_3) = \frac{1}{z^2}+
\underset{\begin{subarray}{c} 
\mathrm{poles}\, \Omega  \\ 
\Omega\neq 0 
\end{subarray}}{\sum}
\left(\frac{1}{(z-\Omega)^2}-\frac{1}{\Omega^2}\right)  
\ee
(cf. \cite{steinmetz2}, Theorem 7.1). The main difference 
with the elliptic case is that the poles of $u$ no longer lie 
on a lattice, and suitable analogues of elliptic integrals to determine the 
positions of the poles are unavailable. 

Nevertheless, %several observations hold. 
the formula 
(\ref{cexp}) for the coefficients in (\ref{series}) still holds in general, 
and as a consequence, from (\ref{orec}), $c_3=0=c_7$ gives 
$\sum_{\Omega\neq 0} \Omega^{-3}=0=\sum_{\Omega\neq 0} \Omega^{-7}$; 
while generically one expects that $c_n\neq 0$ for $n\geq 4$, $n\neq 7$ 
(and similarly for the corresponding power sums). Following the 
elliptic case, by picking a non-zero pole $\Omega_{*}$ such that 
$|\Omega_{*}|$ is minimal, one may write 
\be \label{cnfo} 
c_n = \frac{n-1}{\Omega_{*}^n}\, F_n, \qquad 
F_n = \sum_{\hat{\Omega}\neq 0}\hat{\Omega}^{-n}, 
\ee 
where the sum is over the non-zero poles $\hat{\Omega}=\Omega /\Omega_{*}$ 
of the rescaled solution 
$
u(z;\hat{g}_2,\hat\lm ,\hat{g}_3)$ with % \qquad \mathrm{with} \qquad  
$\hat{g}_2=\Omega_{*}^4g_2$, %\quad 
$\hat\lm = \Omega_{*}^5\lm$, %\quad 
$\hat{g}_3=\Omega_{*}^6g_3$. 
%$$ 
Generically, for the rescaled solution, there should be only one pole 
with modulus 1 (corresponding to $\Omega_{*}$), with all other 
non-zero poles lying outside the unit circle. Assuming this to be the case, 
it follows that 
$\lim_{n \to\infty}F_n = 1$, whence %\qquad \mathrm{whence} \qquad 
$\Omega_{*}=  \lim_{n \to\infty}\frac{c_n}{c_{n+1}}$.  
Thus, in the generic situation, %we find that 
numerical iteration of (\ref{rec1}) allows the pole nearest to $z=0$ to be calculated very efficiently.  
 
The question of precisely when a non-generic solution 
can arise, with two or more poles $\Omega$ 
having the same modulus $|\Omega_{*}|$, 
seems to be a very difficult one. However, there is one situation 
where we know this to be so, namely when $g_2=0=g_3$: we call this 
the {\it pentagonal case}. It corresponds to a solution which 
is invariant under the order 5 scaling symmetry of $P_I$, 
which is generated by taking $z\to e^{2\pi\i' /5 }z$, 
$u\to e^{6\pi\i' /5 }u$ in equation (\ref{p1orig}). 
There are precisely two such solutions, and it appears that so far 
these are the only solutions of $P_I$ for which the 
monodromy data can be explicitly related to the initial conditions 
at $z=0$ \cite{kitaev}:  one has 
a double pole, and the other has a triple zero at the origin; it is the former 
symmetric solution which is of interest to us here. 

\noindent {\bf Pentagonal case:} When $g_2$ and $g_3$ both vanish,  
it follows from (\ref{orec}) that $c_n=0$ unless 
$n=0\,\bmod \, 5$, so that the Laurent series (\ref{series}) 
for %the solution 
$\tl{u}(z)=u(z;0,\lm ,0)$  
takes the form 
\begin{equation}\label{sol1} 
%u^p(x)
\tilde{u} (z) 
=\frac{1}{z^2}+ 
%\frac{5}{6} 
\sum_{n=1}\frac{5n-1}{\gamma^n}\tilde{F}_n 
\lm^n z^{5n-2}, 
\end{equation}
where we set 
$v_n = c_{5n}=\frac{(5n-1)}{\gamma^n}\, \tilde{F}_n$.  
%Above we introduced %It is convenient to calculate using 
The coefficients ${v}_n$ %which 
satisfy the recurrence 
\begin{equation}\label{recvn} 
v_n=\frac{6}{(5n+1)(5n-6)}\sum_{k=1}^{n-1}v_{k}v_{n-k}, \qquad n\geq 2, %2,3... \\
\qquad \mathrm{with} \quad v_1=1. 
\end{equation} 
In this case, the non-zero 
poles of the solution (\ref{sol1}) lie on regular pentagons centred 
at $z=0$: when $\Omega$ is a pole of $\tl{u}$, then so 
is $e^{\frac{2\pi\i' j}{5}}\Omega$ for $j=1,2,3,4$. 
We have $F_n=0$ if $n\neq 0 \,\bmod\,5$,  and from (\ref{cnfo}), with 
5 poles on the unit circle,  
$\lim_{n\to\infty} F_{5n} =5$.  
%twenty decimal places at $n=27$, whether the corresponding value  is:
%Calculating up to 
At $n=30$ we find a value of $\gm$ apparently correct to 23 decimal places:
\be\label{gam}
\gm = \lim_{n\to\infty}\frac{(5n+4){v}_n}{(5n-1){v}_{n+1}} 
\approx 18.32138268472483887119960... .
\ee 
%appears to be  
Given this accuracy in %the value of 
$\gm$, the %values of the 
quantities $\tl{F}_n = F_{5n}$ can then be calculated, 
and should converge exponentially fast to the value 5; see Table \ref{tab3}. 
The value $\gm^{1/5}\approx 1.788923$ gives the non-zero real pole 
of the solution $u(z;0,1,0)$ closest to the origin, and the radius of 
convergence of the series (\ref{sol1}); for an independent verification of this 
%result is given in 
see section 4. 

The preceding results on limiting values of the coefficients in the symmetric cases %{MLpi} {cnfo}
can be summarized as follows. %when 
\begin{propn} 
The solutions (\ref{MLpi}) for which only one of the parameters $g_2,\lambda ,g_3$ 
is non-zero are invariant under the scaling $z\to \xi z$, $u\to \xi^{-2}u$, where 
$\xi$ is a $k$th root of unity for $k=4,5,6$ respectively, and the normalized 
coefficients (\ref{cnfo}) satisfy 
$$ 
\lim_{n\to\infty}F_{kn}=k. 
$$
\end{propn} 
\begin{prf}
The normalized series $F_n=\sum_{\hat{\Omega}\neq 0}\hat{\Omega}^{-n}$ are absolutely convergent 
(uniformly in $n$), and for a solution invariant under the symmetry of order $k$, 
the sum of the series vanishes unless $k|n$. With the chosen normalization, there are  $k$ poles 
at the roots of unity $\hat{\Omega}=\xi^j$, $j=0,\ldots ,k-1$. Given that there are 
no other poles on the unit circle, we find  $F_{kn} = k + \sum_{|\hat{\Omega}|>1}\hat{\Omega}^{-kn}$, 
and hence (upon selecting the term $\hat{\Omega}=\tilde{\Omega}$ with smallest modulus outside the unit circle) 
$
|F_{kn}-k| = \left| \tilde{\Omega}^{-kn} \sum_{|\hat{\Omega}|>1}\left( \tilde{\Omega}/\hat{\Omega}\right)^{kn}\right| 
\leq |\tilde{\Omega}|^{-kn}\sum_{|\hat{\Omega}|>1}\left| \tilde{\Omega}/\hat{\Omega}\right|^{k}\to 0
$
exponentially fast as $n\to\infty$. 
\end{prf} 
The only statement above that requires further justification 
is the assertion that there are no other poles on the unit circle for the case $k=5$ 
(Note that  for $k=4,6$ this is obvious from the properties of the pole lattice.) A numerical 
verification of this fact %, obtained by truncating the Taylor series of the tau-function to calculate its zeros, 
is given in the next section 
(see Figure 1). 
We conclude this section with some observations on the cases $g_2=0$ 
and $g_3=0$. 
\begin{table}
\begin{center}
\begin{tabular}{|c|c|}
\hline
$n$ & $\tl{F}_n$\\
\hline \hline
1 & 4.58034567118120971779 \\ %5.49641480541745166135 \\
2 & 5.08595550727477491732 \\ %6.10314660872972990079 \\
3 & 4.99187877676419618477 \\ %5.99025453211703542173 \\
4 & 5.00112762186482314743 \\ %6.00135314623778777691 \\
5 & 4.99986996982708054870 \\ %5.99984396379249665844 \\
6 & 5.00001616272241466829 \\ %6.00001939526689760195 \\
... & ...\\
11 & 4.99999999957591996469 \\ %5.99999999949110395762 \\
12 & 5.00000000005151463070 \\ %6.00000000006181755684 \\
13 & 4.99999999999374379484 \\ %5.99999999999249255381 \\
14 & 5.00000000000075986460 \\ %6.00000000000091183753 \\
%...& ...\\
\hline
\end{tabular}
\caption{The oscillating values of $\tl{F}_n=F_{5n}$.}
\label{tab3}
\end{center}
\end{table}

\subsubsection{The case $g_2=0$}\label{subsec2} 

%c_4=\frac{g_2}{20}\doteq g, \qquad c_6\doteq a

To consider the form of the quantities $c_n$ when $g_2=0$, it is convenient 
to use the parameter $\alpha =g_3/28$. Thus we have the recursion 
(\ref{rec1}) with $c_0=1$, $c_1=c_2=c_3=c_4=0$, $c_5=\lm$ and $c_6=\alpha $, and find that  
the structure of the iterates depends on the index mod 5, so that  
%It is also quite straightforward to insert $\lm$ into the formulae, the result for the term $c_k$ being:
\begin{equation}
c_{5n+p}=\sum_{m=0}^{\left[\frac{n-p}{6}\right]}c^{(m)}_{5n+p}\, \alpha^{5m+p}\lm^{n-6m-p} , \qquad p=0,\ldots,4 , 
\end{equation}
for certain rational numbers $c_{5n+p}^{(m)}$, which satisfy quadratic recurrence relations.

To begin with, we examine  $c^{(0)}_{5n}$, which are the coefficients of $\alpha^0$. Upon setting 
$v_n\doteq  c_{5n}^{(0)}$, we see that $v_n$ satisfies (\ref{recvn}). %corresponding via the scaling $\tl{v}_n = \frac{6}{25}v_n$. 
The next terms to consider are $c^{(0)}_{5n+1}$, i.e. the coefficients of $\alpha^1$, and upon setting 
$w_{n}\doteq c^{(0)}_{5n+1}$ we find  
\begin{equation}\label{recwn} %\left\{\begin{aligned}
%& 
w_n=\frac{12}{(5n+2)(5n-5)}\sum_{k=1}^{n-1}w_{k}v_{n-k}, \qquad n\geq 2, %=2,3... \\
%& 
\qquad \mathrm{with} \quad w_1=1 .
%\end{aligned}\right.
\end{equation}
Observe that, supposing the sequence of $v_n$ to be known, the above recurrence 
is \emph{linear} in the unknowns (unlike the recurrence (\ref{recvn}) for the $v_n$ themselves). 
This means that the generating function for $w_n$ satisfies a second order linear ODE, within which the 
generating function of the $v_n$ appears as a coefficient.  
 
As we shall now see,  this property extends to the recursive generation of all the sequences $c^{(m)}_{5n+p}$ when $(m,p)\neq(0,0)$. 
The linear recurrence solved by the general term $c^{(m)}_{5n+p}$ can be written explicitly as 
\be\label{reccgen}
c^{(m)}_{5n+p}=K_{n,p}\sum_{k=1}^{n-1}\left(\sum_{j=0}^m\sum_{\ell =0}^p c^{(j)}_{5k+\ell}c^{(m-j)}_{5(n-k)+p-\ell}+\sum_{j=0}^{m-1}\sum_{\ell=1}^{4-p}c^{(j)}_{5k+\ell +p}c^{(m-j-1)}_{5(n-k)-\ell}\right), \quad n\geq 2, 
\ee
where $K_{n,p}$ is given by 
\be\label{Vnp}
K_{n,p}=\frac{6}{(5n+p+1)(5n+p-6)}, \qquad p=0,1,2,3,4. 
\ee
%and again we recall that the values of $p$ are limited to the set $(0,1,2,3,4)$. 
The initial conditions for the equations (\ref{reccgen}) are that $c^{(m)}_k =0$ 
for $1\leq k\leq 9$ and for all $m\geq 0$, {\it except}  
for  $c^{(0)}_5=c^{(0)}_6=1$.  
From (\ref{reccgen}) it appears that 
all the nonlinearity in the problem of determining these coefficients is moved into finding the solution of (\ref{recvn}). 

For the recurrence (\ref{recvn}), note that it is sufficient to determine a single non-vanishing solution, since all other 
solutions can be obtained by the rescaling $v_n\to A^n v_n$, with $A$ arbitrary.  
The $v_n$ correspond precisely to the coefficients in the Laurent expansion (\ref{sol1}) of the 
solution $\tl{u}(z)$ with pentagonal symmetry that was discussed previously. To see an example of how this is related to the other 
coefficients, we consider $w_n = c^{(0)}_{5n+1}$ once more, and introduce the generating function $G(x)= \sum_{n=1}^{\infty}w_nx^{n-1}$. 
From (\ref{recwn}) it follows that $G$ satisfies %the ODE 
$ x\, G''+\frac{12}{5}\, G' =\frac{12}{25}\, G\psi$, 
where $\psi(x)=\sum_{n=1}^{\infty}v_nx^{n-1}$. In fact,  up to some scaling and shifting, $\psi(x)$ is just 
given by the pentagonal solution $\tl{u}$; to be precise,  
$x\psi(x)=\left(x/\lm \right)^{\frac{2}{5}}\tl{u} %\left( 
(\left(x/\lm\right)^{\frac{1}{5}}
)%\right)
-1$. 
In a similar way, via (\ref{reccgen}), 
the generating functions for the other terms $c^{(m)}_{5n+p}$ are related to each other: for example, 
that of $c^{(0)}_{5n+2}$ is related to the generating functions of the sequences $c^{(0)}_{5n}$ and $c^{(0)}_{5n+1}$; 
and that of $c^{(1)}_{5n}$ is related to the generating functions of the sequences $c^{(0)}_{5n+p}$, $p=0,\ldots ,4$, and so on.

\subsubsection{Case $g_3=0$ }\label{subsec3}
For the case $g_3=0$, it is helpful to introduce the parameter $\beta =g_2/20$, and 
take the solution of (\ref{rec1}) with 
$c_0=1$, $c_1=c_2=c_3=c_6=0$, $c_5=\lm$ and $c_4=\beta$. The coefficients 
of (\ref{series}) now have the structure  
\begin{equation}
c_{5n-p}=\sum_{m=0}^{\left[\frac{n-p}{4}\right]} \hat{c}^{(m)}_{5n-p} \, \beta^{5m+p}\lm^{n-4m-p}, \qquad p=0,\ldots,4,
\end{equation}
for some rational numbers $\hat{c}^{(m)}_{5n-p}$. 

The quantities  $\hat{c}^{(0)}_{5n}$, which appear as coefficients of $\beta^0$, are precisely the same as 
the numbers $v_n=c^{(0)}_{5n}$ found previously, since setting $\beta=0$ just gives the pentagonal 
solution $\tl{u}$ with the expansion (\ref{sol1}). The next terms to consider are $\hat{w}_n\doteq \hat{c}^{(0)}_{5n-1}$, 
the coefficients of $\beta^1$, 
which satisfy the recurrence
\begin{equation}\label{recxin} %\left\{\begin{aligned}
%& 
\hat{w}_n=\frac{12}{(5n)(5n-7)}\sum_{k=1}^{n-1}\hat{w}_{k}{v}_{n-k} \qquad n\geq 2, %3... \\
%& 
\qquad \mathrm{with} \quad \hat{w}_1=1. 
%\end{aligned}\right.
\end{equation}
Just as for the case $g_2=0$, once $v_n$ is given, this recurrence is {\it linear} in the unknowns 
$\hat{w}_k$. 
Analogously to equation (\ref{reccgen}), it is possible to write down the general recurrence solved by the term $\hat{c}^{(m)}_{5n-p}$, 
the result being 
\be\label{recbgen}
\hat{c}^{(m)}_{5n-p}=K_{n,-p}\left(\sum_{j=0}^m\sum_{\ell =0}^p\sum_{k=1}^{n-1}\hat{c}^{(j)}_{5k-\ell}\hat{c}^{(m-j)}_{5(n-k)+\ell -p} 
+\sum_{j=0}^{m-1}\sum_{\ell =1}^{4-p}\sum_{k=1}^{n}\hat{c}^{(j)}_{5k-\ell -p}\hat{c}^{(m-j)}_{5(n-k)+\ell }\right) , 
\ee
where again $K_{n,p}$ is given by (\ref{Vnp}), but above it appears with $p\to -p$ compared with 
(\ref{reccgen}). 
The formula (\ref{recbgen}) holds for 
%%%%$n\geq 3$ when $p=4$, and for 
$n\geq 2$,  %%%%in all other cases, 
and the initial conditions are given by $\hat{c}^{(m)}_k=0$ for $1\leq k \leq 6$ and for all $m\geq 0$, {\it except} 
%b^z_2=...=b^z_6=0$ 
for $\hat{c}^{(0)}_4=\hat{c}^{(0)}_5=1$. 
Similarly to the situation for $g_2=0$, linear ODEs for the generating functions of 
the rational numbers $\hat{c}^{(m)}_{5n-p}$ for $(m,p)\neq (0,0)$ can be constructed recursively, 
once the numbers $v_n$ are known.

\section{Expansion of the tau-function} 

The tau-function for (\ref{PI}), related to $u$ by (\ref{tau}),  
satisfies a fourth order differential equation  which is homogeneous of degree two. %and 
%can be 
It is written in Hirota bilinear form as %equation 
\be \label{bil} 
D_z^4 \uptau\cdot \uptau -(12\lm z + g_2)\uptau^2 =0, 
\ee 
where the Hirota derivative $D_z$ is  defined by 
$ D_z^n f\cdot g (z) = \left( \frac{d}{dz}-\frac{d}{dz'}\right)^n f(z)g(z') |_{z'=z}$.  

Painlev\'e analysis can be applied directly to the equation (\ref{bil}), 
%making 
expanding around a simple zero at $z=0$, 
corresponding to a double pole in $u$ there. Seeking resonances by taking 
$
\uptau \sim z + \epsilon z^r 
$
yields $r=-1,0,1,6$: the value $-1$ is the movable position of the singularity, as usual, 
and the values $r=0,1$ correspond to the two free parameters $a,b$ for the gauge transformations 
(\ref{gauge}), which  leave the equation (\ref{bil}) invariant; this leaves only $r=6$, which 
is equivalent to the freedom to choose $c_6$ (or the parameter $g_3$) in the Laurent series (\ref{series}).  

We would like to determine the Taylor series of the tau-function around 
$z=0$, with $\uptau (0)=0$, 
since  (due to the fact that $\uptau$ is an entire function), this provides a global representation of 
the solution of (\ref{PI}), via the formula (\ref{tau}).  
If we fix the gauge, we can always choose the coefficient of $z$ to be 1, and set the coefficient of 
$z^2$ to be $0$, which results in a series of the form 
\be\label{taus} 
\uptau (z) = z + \sum_{n=2}^\infty C_{n} z^{n+1}. 
\ee 
This is the analogue of the power series for the Weierstrass sigma function $\sigma (z)$, which satisfies 
the bilinear equation (\ref{bil}) with $\lm =0$. Using the same method as in \cite{ee}, with %and making use of 
$$ 
D_z^4 z^j\cdot z^k =b_{j,k}z^{j+k-4}, \qquad \mathrm{where} 
\qquad 
b_{j,k}=4!\sum_{\ell =0}^4 (-1)^\ell \left(\begin{array}{c} j \\ \ell \end{array}\right) \, 
\left(\begin{array}{c} k \\ 4 - \ell \end{array}\right) , 
$$ 
it is 
straightforward to establish the following.

\begin{thm}\label{cthm} 
The coefficients in the expansion (\ref{taus}) of the tau-function belong to 
$\mathbb{Q}[g_2,\lm ,g_3]$, being uniquely 
determined by the recursion 
$$ 
n(n^2-1)(n-6)C_n =-\frac{1}{2}\sum_{j=1}^{n-1}b_{j+1,n-j+1}C_jC_{n-j} 
+ \frac{1}{2}g_2\sum_{j=0}^{n-4}C_jC_{n-4-j}+6\lm  \sum_{j=0}^{n-5}C_jC_{n-5-j}
$$ 
subject to fixing $C_0=1,C_1=0,C_6=-g_3/840$.  
Each %coefficient 
$C_n$ is a weighted homogeneous polynomial of total degree $n$ %where 
in the arguments $g_2,\lm ,g_3$ with weights $4,5,6$ respectively. 
\end{thm}    
\begin{rem} The choice of gauge is specific to the zero at $z=0$. 
If we let  $\uptau (z; g_2,\lm ,g_3 )$ denote the function 
given by (\ref{taus}), 
and then expand around another zero at $z=\Omega \neq 0$, we obtain 
the formula 
$$ 
\uptau (z; g_2,\lm ,g_3 ) = A e^{B (z-\Omega)}  
\uptau (z-\Omega; g_2+12\lm\Omega ,\lm ,\hat{g}_3 )   , 
$$  
where in principle $\hat{g}_3$, $A=\uptau '(\Omega )$ and 
$B=\frac{1}{2}\uptau ''(\Omega )/\uptau '(\Omega )$ 
all depend on $g_2,\lm ,g_3$ and $\Omega$. The quasiperiodicity 
of $\sigma (z)$ % the sigma function 
under shifting by a period is a special case of this. %the latter formula.  
\end{rem} 
 
For the case of the sigma function in \cite{ee} 
it is noted that the value of the coefficient $C_6$ must be given
appropriately  in terms of $g_3$  
in order to be consistent with the Laurent series for $\wp$ satisfying 
(\ref{wpeq}). The same holds for $\lm \neq 0$, upon requiring  
consistency of (\ref{taus}) with (\ref{series}), as is seen 
by noting that the expression (\ref{ham}) for the Hamiltonian 
gives a first integral for the bilinear equation (\ref{bil}), taking  
$v=u'$ with $u$ given by (\ref{tau}) and $h$ given by (\ref{htau}). 
Then the formula (\ref{ham}) can be rewritten as a third order equation for $\uptau$ 
that is homogeneous of degree four, namely 
\be \label{quad} 
\begin{array}{lc} 
\uptau^2(\uptau ''')^2 -6\uptau\uptau '\uptau ''\uptau '''  
+4(\uptau ')^3\uptau ''' +4\uptau (\uptau '')^3 -3(\uptau ' \uptau '')^2 
&\\ 
-g_2\, \uptau^2 \Big(\uptau\uptau '' -(\uptau')^2\Big) 
-12\lm   \Big(z(\uptau^3\uptau '' -(\uptau\uptau')^2) -\uptau^3\uptau ' \Big) 
+ g_3\, \uptau^4 & = 0.  
\end{array} 
\ee  
This equation immediately yields another recurrence for the coefficients 
of %the Taylor expansion for 
$\uptau$. %(\ref{taus}).  
\begin{thm}\label{cthm2} 
Subject to fixing $C_0=1,C_1=0$, the coefficients in the 
expansion (\ref{taus}) are 
uniquely determined by the recursion 
$$ 
n(n^2-1)\, C_n =-\underset{\begin{subarray}{c} 
j+k+\ell+m=n \\ 
1\leq j,k,\ell ,m \leq n-1 
\end{subarray}}{\sum}\tilde{b}_{j,k,\ell ,m}\, C_jC_{k}C_\ell C_m  
+ g_2\sum^{(4)} +12\lm \sum^{(5)} -g_3 \sum^{(6)}, 
$$ 
where 
$$ 
\tilde{b}_{j,k,\ell ,m} =j(j+1)(k+1)\Big( 
(j-1)\big( k(k-6\ell -7) +4(\ell +1)(m+1) \big) +k(\ell +1)(4\ell -3m-3)\Big)  
,
$$ 
and 
$$ 
\sum^{(4)}=\underset{\begin{subarray}{c} 
j+k+\ell+m=n-4 \\ 
0\leq j,k,\ell ,m \leq n-4 
\end{subarray}}{\sum} (j+1)(j-k-1)C_jC_{k}C_\ell C_m , 
$$ 
$$ 
\sum^{(5)}=\underset{\begin{subarray}{c} 
j+k+\ell+m=n-5 \\ 
0\leq j,k,\ell ,m \leq n-5 
\end{subarray}}{\sum} (j+1)(j-k-2)\, C_jC_{k}C_\ell C_m , 
\quad \sum^{(6)}=\underset{\begin{subarray}{c} 
j+k+\ell+m=n-6 \\ 
0\leq j,k,\ell ,m \leq n-6 
\end{subarray}}{\sum}  C_jC_{k}C_\ell C_m .  
$$ 
\end{thm} 
\begin{figure}\centering
\scalebox{0.4}[0.55]{
\includegraphics[angle=270]{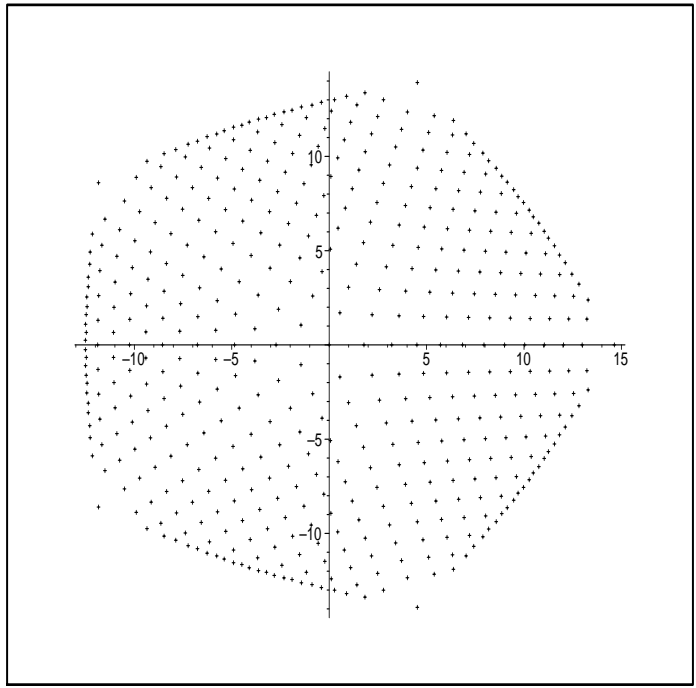}}
\caption{\small{Approximation to the poles of the solution of (\ref{PI}) with $g_2=g_3=0$.}}
\label{pentangle} 
\end{figure}
Both of the recursions for $C_n$ given above can be iterated very rapidly 
to produce polynomial approximations to the tau-function, for any choice 
of the values $g_2,\lm , g_3$. As an example, to generate 
Figure \ref{pentangle} we took the values $g_2=0=g_3$ and $\lm =1$, 
corresponding to the solution with pentagonal symmetry that was 
considered in the previous section, and calculated the first 101 non-zero 
terms in the expansion (\ref{taus}), which contains only terms 
of the form $z^{5n+1}$, beginning with  
$
\uptau (z) = z -\frac{z^6}{20}-\frac{7z^{11}}{26400} 
+ \frac{z^{16}}{1232000}+ \frac{83z^{21}}{117976320000}+\ldots $.  
In the figure we have plotted the roots of the truncated series, 
which is of the form $P_{501}(z)=z\, \hat{P}_{100} (w)$, where 
$\hat{P}_{100}$ is a polynomial of degree 100 in $w=z^5$. The roots of 
the truncated series approximate the zeros of $\uptau$ with 
increasing accuracy as more terms are added. The smallest root 
of $\hat{P}_{100}$ is a real number $\gamma$, agreeing with 
the digits of the numerical value (\ref{gam}) to the same 
accuracy as before. In fact, it appears easier to get improved 
approximations to this value by truncating the Taylor series 
than by calculating the ratios of the coefficients $c_n$ 
in the Laurent expansion.

Before concluding this section, we present another representation for the 
series (\ref{taus}), which displays some interesting arithmetical 
features of the coefficients. Note that, since the $C_n$ are weighted 
homogeneous, the tau-function can be written in the form of a triple sum
\be\label{taua} 
\uptau(z) = \sum_{\ell , m,n\geq 0}A_{\ell ,m,n} \, (\frac{1}{2}g_2)^\ell 
(6\lm )^m (2g_3)^n \frac{z^{4\ell + 5m + 6n+1}}{(4\ell + 5m + 6n+1)!},  
\ee 
for certain rational numbers $A_{\ell ,m,n}$. 
The above formula is motivated by the result of Weierstrass 
\cite{weierstrass}, who showed 
that the sigma function can be expressed as 
$
\sigma (z) = \sum_{m,n\geq 0}a_{m,n} \, (\frac{1}{2}g_2)^m 
 (2g_3)^n \frac{z^{4m + 6n+1}}{(4m + 6n+1)!} $. 
From comparison of the two series it is clear  that upon setting $\lm =0$ we have 
$ a_{m,n} = A_{m,0,n}$ for $m,n\geq 0$.  
Weierstrass used (\ref{wp}) and (\ref{wpeq}), together with various 
modular relations, to show that the sigma function satisfies two linear 
partial differential equations, from which he obtained 
a linear recurrence relation for the coefficients $a_{m,n}$. 
The first of these PDEs just follows 
from Euler's theorem on homogeneous functions, and thus extends 
to the tau-function when $\lm\neq 0$ also: 
$$ 
\left( 4g_2\frac{\partial}{\partial g_2} 
+5\lm \frac{\partial}{\partial \lm} + 6g_3\frac{\partial}{\partial g_3} 
-z\frac{\partial}{\partial z} +1 \right) \, \uptau = 0.
$$ 
(Note that $z$ has weight $-1$, so overall $\uptau$ has the same weight.)
However, the other PDE that Weierstrass found for $\sigma (z)$ is of 
second order in $z$ (containing also first derivatives with respect 
to $g_2,g_3$);  it is equivalent to the fact that the corresponding 
elliptic theta function satisfies the heat equation. 
For the tau-function of $P_I$ we do not expect any linear PDE of this kind. 
Nevertheless, 
the Hirota bilinear equation (\ref{bil}) provides recursive relations
for the coefficients $A_{\ell ,m,n}$. 
\begin{thm}
 Subject to fixing $A_{0,0,0}=1$, $A_{0,0,1}=-3$, the coefficients in (\ref{taua}) 
are completely determined by the recursion 
$$ 
\begin{array}{ccl} 
\frac{s(s-1)(s-2)(s-7)}{s!}\, A_{\ell , m,n} & = & -\frac{1}{2}  \underset{\begin{subarray}{c} 
\ell_1+\ell_2=\ell \\ 
m_1+m_2=m \\
n_1+n_2=n 
\end{subarray}}{\sum '} \frac{b_{s_1,s_2}}{s_1!s_2!} A_{\ell_1 ,m_1,n_1}A_{\ell_2 ,m_2,n_2} 
\\ 
&& +  %%%% 6 %%%% FACTOR OF 6 removed with new scaling!!!
\underset{\begin{subarray}{c} 
\ell_1+\ell_2=\ell \\ 
m_1+m_2=m-1 \\
n_1+n_2=n 
\end{subarray}}{\sum } \frac{ A_{\ell_1 ,m_1,n_1}A_{\ell_2 ,m_2,n_2} }{s_1!s_2!}
+ 
\underset{\begin{subarray}{c} 
\ell_1+\ell_2=\ell -1\\ 
m_1+m_2=m \\
n_1+n_2=n 
\end{subarray}}{\sum } \frac{ A_{\ell_1 ,m_1,n_1}A_{\ell_2 ,m_2,n_2} }{s_1!s_2!}
,
\end{array}  
$$
where $s=4\ell +5m + 6n+1$, $s_j=4\ell_j +5m_j + 6n_j+1$ for $j=1,2$, 
and $\sum '$ denotes that the terms with $(\ell_1 ,m_1,n_1) = (0,0,0)$ or $(\ell ,m,n)$ are omitted 
from the sum. 
\end{thm} 
\begin{rem}Another recursion for $A_{\ell ,m,n}$ can be obtained from (\ref{quad}). 
\end{rem} 
The first few coefficients are given below in the form of $3\times 3$ 
matrices ${\bf M}^{(m)}$ 
whose $(j,k)$ entry is $A_{j-1,m,k-1}$ for $m=0,1,2$: 
$$ 
{\bf M}^{(0)}=\left(\begin{array}{ccc} 1 & -3 & -54 \\ -1 & -18 & 4968 \\ -9 & 513 & 257580 \end{array}\right) , 
 \quad 
{\bf M}^{(1)}=\left(\begin{array}{ccc} 
%NEW COEFFS SCALED BY 6 
-6 &  -216 &  89424 \\ -84 &  18720 &  5786640 \\ 1650 &  1358640 &  1168920720   
%-36 & -1296 & 536544 \\ -504 &  112320 &  34719840 \\  9900 &  8151840 &  7013524320 
\end{array}\right) , 
$$ 
$$ 
{\bf M}^{(2)}=\left(\begin{array}{ccc} 
%NEW COEFFS SCALED BY 36 
-294 &  144144 & 47585880 \\ 18774 & 15053040 &  22914336240 \\ 1112436 &  3160803600 & -2734614623160  
%-10584 & 5189184 & 1713091680 \\ 675864 & 541909440 & 824916104640 \\ 40047696 & 113788929600 &  %-98446126433760  
\end{array}\right) . 
$$

Due to recent 
results of Onishi \cite{onishi}, it is known that the coefficients $a_{m,n}$ in the 
expansion of the sigma function are all integers, and analogous results 
have been proved for sigma functions of some higher genus curves 
\cite{eilbeck}. This suggests  
\begin{conje} 
The series (\ref{taua}) has 
$
A_{\ell ,m,n} \in \mathbb{Z} \qquad \forall \ell ,m,n \geq 0. 
$ 
\end{conje} 
We have verified this conjecture for the first few hundred coefficients. 

\section{Conclusions}  
The properties of the exact series expansions for the solution of $P_I$ 
and its tau-function are tantalizingly close to those of the 
analogous Weierstrass functions. 
As well as their potential uses for numerical calculations,  it would be 
interesting to use the recursions in Theorems \ref{cthm} 
and \ref{cthm2} to prove directly that $\uptau$ is holomorphic 
(hence providing yet another proof of the fact that all solutions 
of $P_I$ are meromorphic). 
Finally, we should mention that all of the Painlev\'e equations 
give deformations of elliptic functions, and they all have   
associated tau-functions and bilinear equations \cite{hietarinta}, so 
the same methods could be applied to $P_{II-VI}$. 
We also note that 
exact series formulae for $P_{VI}$ tau-functions have been 
obtained recently from a completely different perspective, 
in terms of correlation functions in 
2D conformal field theory \cite{lisov}.  

\noindent {\bf Acknowledgments.} ANWH is grateful to Martin Kruskal, who 
suggested several years ago that it might help to consider the Taylor series 
for the tau-function. He also thanks Sasha Veselov for discussions at the 
MISGAM workshop in Berlin in 2005, concerning a similar suggestion 
due to Alexei Shabat. FZ wishes 
to acknowledge the financial support of the Istituto Nazionale di Alta 
Matematica, in the form of %who provided him with 
an INdAM-COFUND Marie Curie fellowship. Both authors are grateful to the anonymous 
referee who helped to improve the manuscript.  
 
%\tiny 

\end{document}